\documentclass[reqno,centertags]{amsart}
\usepackage{pstricks-add,pstricks}
\usepackage{amsmath,amsfonts,amssymb,amsxtra,latexsym,amscd,enumerate,amsthm,hyperref}
\usepackage{dsfont}
\usepackage{latexsym}
\usepackage{epsfig}
\usepackage{comment, color}
\usepackage[final]{showkeys}
\usepackage{cleveref}

\newtheorem{d1}{Definition}
\newtheorem{t1}{Theorem}
\newtheorem{l1}{Lemma}
\newtheorem{p1}{Proposition}

\newtheorem*{r1}{Remark}

\newcommand{\eps}{\epsilon}
\newcommand{\jb}[1]{\langle #1 \rangle}

\newcommand{\ind}[1]{\mathbf{1}_{#1}}
\newcommand{\R}{\mathbb{R}} \newcommand{\Z}{\mathbb{Z}}
 \newcommand{\C}{\mathbb{C}}
\newcommand{\N}{\mathbb{N}}
\newcommand{\F}{\mathcal{F}}
\newcommand{\supp}{\operatorname{supp}}

\newcommand{\lb}{\langle}
\newcommand{\rb}{\rangle}
\newcommand{\ls}{\lesssim}

\specialcomment{note}{\begingroup\sffamily\color{blue} }{\endgroup}

\crefname{t1}{Theorem}{Theorems}
\crefname{p1}{Proposition}{Propositions}
\crefname{l1}{Lemma}{Lemmas}

\begin{document}

\title{Notes on Delort's KG paper}
\author{Tobias~Schottdorf} \email{schottdo@math.uni-bonn.de}

\address{
  Rheinische Friedrich-Wilhelms-Universit\"at Bonn, Mathematisches Institut,
  Endenicher Allee 60, 53115 Bonn, Germany.}

\title[quadratic Klein-Gordon equations]{
  Global existence without decay for quadratic Klein-Gordon equations}

\begin{abstract}
  The Cauchy problem for quadratic Klein-Gordon systems \[ (\Box + m_i^2)u_i = N(u_1, \ldots, u_K) \qquad i=1,\ldots,K \] is considered in two spatial dimensions and higher under a suitable non-resonance condition on the masses $\{ m_i \}$, including the main case of equal masses. A global well-posedness and scattering result is proven for small $H^s$ data, $s \geq \max(\frac{1}{2}, \frac{n-2}{2})$, using the contraction mapping technique in $U^2/V^2$ based spaces. The result demonstrates that global results and scattering hold in this setting without the need to impose strong decay on the initial data.
\end{abstract}



\maketitle

\setcounter{tocdepth}{1}
\tableofcontents

\section{Introduction and main results}
From the late 1970s on, there has been a lot of progress on questions of global existence and blow-up for equations and systems of the type
\begin{equation}\label{nl}
\begin{aligned} (\Box + m^2)u(t,x)  &= F_p(u(t,x)),\qquad (t, x) \in [0, T) \times \R^n \\ u(0, x)&=f(x)\\ \partial_t u(0,x) &= g(x)
\end{aligned}
\end{equation}
where $m \geq 0$, $\Box = \partial_{tt} - \Delta$, $u$ is scalar or vector-valued, and $F_p$ is a power-type nonlinearity of order $p > 0$, i.e. $|\partial^j F_p(s)| \sim |s|^{p-j}$ $(j \leq p)$ together with a similar condition for differences.\\

An optimistic energy heuristic based on the decay of free solutions leads to a first guess that global existence from small data could hold for 
\[ \begin{aligned}p &> 1 + \frac{2}{n}\qquad &\text{ if }m>0 \\ p &> 1 + \frac{2}{n-1}\qquad &\text{ if }m=0. \\\end{aligned}\]

We shall be interested primarily in the first case $m > 0$ but summarize the massless version $m = 0$ briefly. As it turns out, that case - where one is dealing with a nonlinear wave equation - is somewhat singular in the sense that the above heuristic is incorrect. Instead, the decisive role is played by a larger number commonly known as the \textit{Strauss exponent}, the positive root $\gamma = \gamma(n)$ of
\begin{equation}\label{strauss} \frac{n}{2}\frac{\gamma - 1}{\gamma + 1} = \frac{\gamma}{2}. \end{equation}
Note $\gamma(1) \sim 3.56$, $\gamma(2) \sim 2.41$, $\gamma(3) = 2$, $\gamma(4) = 1.78$, $\gamma(\infty) = 1$ and \[ 1+ \frac{2}{n} < \gamma < 1 + \frac{4}{n}. \]

More precisely, for $m = 0$, $\gamma(n-1)$ is a threshold power such that for \eqref{nl} we have the following dichotomy:
If $p > \gamma(n-1)$, then small, smooth and localized data lead to global solutions. In the other case $p \leq \gamma(n-1)$, one can find such data blowing up in finite time. This conjecture-turned-theorem\footnote{modulo the endpoint, Schaeffer \cite{MR824205} confirmed this for $n=2$ and Glassey \cite{MR631631} subsequently proved finite time blow-up for the critical cases $p = \gamma(n-1)$, in two and three dimensions. For larger $n$, blow-up from small data below $\gamma(n-1)$, was subsequently confirmed by Sideris \cite{MR744303}, while the positive part is due to Zhou \cite{MR1331521}, Lindblad and Sogge \cite{MR1408499}, Georgiev, Lindblad and Sogge \cite{MR1481816}. Finally,  Yordanov and Zhang \cite{MR2195336} proved blow-up also when $p = \gamma(n)$ for the open cases $n \geq 4$} goes back to Strauss, who based his prediction on results by John in 3D \cite{MR535704} and his own work. Hence, for the wave equation case $m = 0$, there is a very clear dichotomy between global solutions and finite time blow-up, indicated by $\gamma(n-1)$.\\

For $m > 0$, where one is dealing with a nonlinear Klein-Gordon equation, the picture is less clear, and in particular the role of $\gamma(n)$. This is somewhat curious since $\gamma(n)$ seems to first have arosen in Strauss' work \cite{MR614228} on scattering in the case $m > 0$. The spaces used in that work are based on the $t^{-\frac{n}{2}}$ time decay of free solutions, and the Strauss exponent $\gamma(n)$ occurs as a natural threshold below which the nonlinearity $|u|^p$ inherits too little decay in time to close the estimates\footnote{the corresponding decay in $L^{p+1}$ is $t^{-d}$, $d = \frac{n}{2}\frac{p-1}{p+1}$. These parameters follow directly from interpolation between the unitary $L^2 \rightarrow L^2$ and the dispersive $t^{-\frac{n}{2}}$ $L^1 \rightarrow L^\infty$ estimates and hence are reasonable to ask for if one expects the solution to behave like a free wave for large times $t$. Chosing to work with $L^{p+1}$ roughly allows solving in $L^\infty(\R, (1+|t|)^{-d} L^{p+1})$ if $1 > d > \frac{1}{p}$. Now the equation $d = \frac{1}{p}$ is exactly \eqref{strauss}.}.\\

From all of the above, it would seem reasonable to expect $\gamma(n)$ to play a the role of a threshold for global existence, scattering, or both\footnote{note that the relevant number is $\gamma(n)$ when $m \neq 0$} when $m > 0$. First insights were again made in three spatial dimensions first, for quadratic nonlinearities by Klainerman \cite{MR803252} and Shatah\cite{MR803256} independently. Noting that $\gamma(3) = 2$, this corresponds to the missing endpoint in Strauss' work \cite{MR614228}.

However, for $n \leq 3$, advances far below the Strauss exponent all the way up to the energy prediction $p > 1 + \frac{2}{n}$ have been obtained by Lindblad and Sogge \cite{MR1408499}; blow-up for $p < 1 + \frac{2}{n}$ in these dimensions is due to Keel and Tao in \cite{MR1738405} at least when $F_p$ is allowed to depend on first derivatives. Additionally, for $n=2$, even in the critical case $p=1+\frac{2}{n}$	global existence is known \cite{MR1400196}.This gives a fairly concise picture in low dimensions; however, Keel and Tao conjecture that in higher dimensions, $1 + \frac{2}{n}$ is not the correct threshold to global existence when $m > 0$.\\

As far as scattering goes in this case, the Strauss exponent also doesn't seem to be a reliable indicator: Recent results by Hayashi and Naumkin \cite{MR2561936}, \cite{MR2373659} show that there is small data scattering in the optimal range\footnote{\label{glassey} \cite{MR0420031}, \cite{MR0330782}} $p > 1 + \frac{2}{n}$ when $n =1, 2$ and for $1 + \frac{4}{n+2} < p < 1 + \frac{4}{n}$ when $n \geq 3$ with small initial data in weighted $H^s$ norms. Both of these results show that scattering results exist below the Strauss exponent in any dimension with only some mild decay on the initial data.

Hence, the Strauss exponent appears to play a less central role, if any, when the masses are positive, at least if the data are sufficiently localized.\\

Since we are interested in systems, we mention only in passing the Hamiltonian theory around $H^1$ data (see, for instance, \cite{MR780103}). Two consequences of this are global existence from small $H^1 \times L^2$ data for $F_p$, $p > 1$, and large data scattering in the energy space for $\gamma(n) < 1 + \frac{4}{n} < p < 1 + \frac{4}{n-1}$. \\

What all of the previous results have in common is that they impose at least some decay on the initial data, the mildest of which being weights or higher $L^p$ norms. We consider it interesting to study the necessity of such conditions, a direction indicated by Delort and Fang in  \cite{MR1789923}. They consider data only in $H^s$ with a reasonable number of derivatives and prove almost global existence for nonlinearities which are quadratic\footnote{and may contain derivatives under the assumption of a null structure} in all dimensions $n \geq 2$. More recently, a global result \cite{georgievstefanov} in this spirit has appeared in the most difficult case $n = 2$ for $H^{1+\varepsilon}$ initial data, which shows that decay is not needed in the endpoint case of the $1 + \frac{2}{n}$ heuristic when $n = 2$.\\

Our contribution is an improvement of this last result in the framework of $U^2$ and $V^2$ spaces (see \cite{MR2526409})  at low regularity. We obtain global existence, scattering and smooth dependence on the initial data for algebraic quadratic nonlinearities in $u$ in dimensions two and higher. It turns out that a certain ``non-resonance'' condition connected to the applicability of the normal forms method allows for a conceptually clear and efficient proof using our setup.

The main result is the following

\newpage
\begin{t1}\label{mainthm}Let $n \geq 2$, $s \geq \max(\frac{1}{2},\frac{n-2}{2})$, $K \in \N$ and let
\[\begin{aligned} N_1, \ldots, N_K &\in \C[x_1, \ldots, x_K] \text{  homogeneous quadratic polynomials, } \\
 m_1, \ldots, m_K &> 0 \text{ such that } 2\min(\{m_j\}) > \max(\{ m_j \} ). \end{aligned}\]
Then there is $\eps > 0$ such that for initial data
\[ (f_i,g_i) \in H^s \times H^{s-1}:\qquad \|(f_i,g_i)\|_{H^s \times H^{s-1}} \leq \eps,\; i = 1, \ldots, K
\] the system
\begin{equation}\label{nlkg}
\begin{aligned}
(\Box + m_i^2)u_i &= N_i(u_1, \ldots, u_K)\qquad i=1,\ldots,K \\
u_i(0) &= f _i\\
\partial_t u_i(0) & = g_i
\end{aligned}
\end{equation}
has a global solution in $C(\R, H^s) \cap C^1(\R, H^{s-1})$. In addition, the solution depends in Lipschitz fashion on $(f,g)$ and scatters asymptotically as $t \rightarrow \pm \infty$.
Furthermore, it is unique in the smaller spaces $X^s([0, \pm \infty))$ introduced in the next section.
\end{t1}

For the sake of clarity, we will first prove the result in the scalar case $K = 1$, where we can assume $m_1 = 1$ and drop the index of $u_i$ and $N_i$. For most arguments it is clear how they carry over to systems; we add the missing pieces in \cref{diffmass}.

\section{Definitions and preliminaries}

\subsection*{Notation}
We denote the spatial Fourier transform by $\hat{\cdot}$ or $\F_x$ and the Fourier transform in time or all variables by $\F_t$ and $\F_{tx}$, respectively.

Frequencies will be denoted by capital letters $M, N, O, H$ and $L$, which we will assume to be dyadic\footnote{we also consider $0$ a dyadic number, corresponding to frequencies smaller than $1$}, that is of the form $2^k$ for $k \in \N$. We will usually not distinguish frequencies below $1$,  numbers $2^k$, $k \in \N$, and we will point out explicitly when a summation is over $2^\Z$ instead.

For dyadic sums, we write $\sum_N a_N :=a_0 + \sum_{n \in \N} a_{2^n}$ and $\sum_{N\geq M} a_N
:=\sum_{n \in \N: 2^n \geq M} a_{2^n}$ for $M > 0$.  Let $\chi \in
C^\infty_0(-2,2)$ an even, non-negative function such that
$\chi(t)=1$ for $|t|\leq 1$. We define $\psi(t):=\chi(t)-\chi(2t)$ and
$\psi_N:=\psi(N^{-1}\cdot)$ for $N > 0$.  Then, $\sum_{N \in 2^\Z}\psi_N(t)=1$ whenever $t \neq 0$, and we define
\begin{equation*}
  \F_t(Q_N u):=\psi_N \F_t(u)
\end{equation*}
and $Q_{\geq M}=\sum_{N \in 2^\Z: N\geq
  M} Q_N$ as well as $Q_{<M}=I-Q_{\geq M}$.

The free evolution $e^{\pm it\jb{D}}$ on $L^2$ is given as a Fourier multiplier
\[ \F_x({e^{\pm it\jb{D}}f}) = e^{\pm it\jb{\xi}} \hat{f} \]

and we define the frequency and modulation projections \begin{align*}
  \widehat{P_N u}(\tau,\xi)&:=\psi_N(\xi)\widehat{u}(\tau,\xi) \text{ for $N \in 2^\N$} \\
  \widehat{P_0 u}(\tau,\xi)&:=\psi_0(\xi)\widehat{u}(\tau,\xi)\\
  \widehat{Q^\pm_M u}(\tau,\xi)&:=\psi_M(\tau\mp\jb{\xi})
  \widehat{u}(\tau,\xi).
\end{align*}
where $\psi_0 = 1 - \sum_{N \geq 1} \psi_N$.

\subsection*{$\mathbf U^2$ and $\mathbf V^2$ spaces}
We review here in condensed form the necessary statements entailing $U^2$ and $V^2$ spaces, taken almost verbatim from \cite{MR2526409}, where detailed proofs can be found\footnote{note also the erratum \cite{MR2629889}}. In what follows, we always have $\frac{1}{p} + \frac{1}{p'} = 1$.\\

Let $\mathcal{Z}$ be the set of finite partitions
$-\infty=t_0<t_1<\ldots<t_K=\infty$ and let $\mathcal{Z}_0$ be the set
of finite partitions $-\infty<t_0<t_1<\ldots<t_K<\infty$.  In the
following, we consider functions taking values in $L^2:=L^2(\R^d;\C)$,
but in the general part of this section $L^2$ may be replaced by an
arbitrary Hilbert space.
\begin{d1}\label{def:u}
  Let $1\leq p <\infty$. For $\{t_k\}_{k=0}^K \in \mathcal{Z}$ and
  $\{\phi_k\}_{k=0}^{K-1} \subset L^2$ with
  $\sum_{k=0}^{K-1}\|\phi_k\|_{L^2}^p=1$ and $\phi_0=0$ we call the
  function $a:\R \to L^2$ given by
  \begin{equation*}
    a=\sum_{k=1}^{K} \mathds{1}_{[t_{k-1},t_{k})}\phi_{k-1}
  \end{equation*}
  a $U^p$-atom. Furthermore, we define the atomic space
  \begin{equation*}
    U^p:=\left\{u = \sum_{j=1}^\infty \lambda_j a_j \text{ in } L^\infty(\R,L^2) \;\Big|\; a_j \text{ $U^p$-atom},\;
      \lambda_j\in \C \text{ s.th. } \sum_{j=1}^\infty |\lambda_j|<\infty \right\}
  \end{equation*}
  with norm
  \begin{equation}\label{eq:norm_u}
    \|u\|_{U^p}:=\inf \left\{\sum_{j=1}^\infty |\lambda_j|
      \;\Big|\; u=\sum_{j=1}^\infty \lambda_j a_j,
      \,\lambda_j\in \C,\; a_j \text{ $U^p$-atom}\right\}.
  \end{equation}
\end{d1}

\begin{p1}\label{prop:u}
  Let $1\leq p< q < \infty$.
  \begin{enumerate}
  \item\label{it:u_banach} $U^p$ is a Banach space.
  \item\label{it:u_emb} The embeddings $U^p\subset U^q\subset
    L^\infty(\R;L^2)$ are continuous.
  \item\label{it:u_right_cont} For $u\in U^p$ it holds
    $\lim_{t\downarrow t_0} \|u(t)-u(t_0)\|_{L^2}=0$, i.e.\ every $u\in
    U^p$ is right-continuous.
  \item\label{it:u_limits} $u(-\infty):=\lim_{t\to - \infty}u(t)=0$,
    $u(\infty):=\lim_{t\to \infty}u(t)$ exists.
  \item\label{it:u_c} The closed subspace $U^p_c$ of all continuous
    functions in $U^p$ is a Banach space.
  \end{enumerate}
\end{p1}

\begin{d1}\label{def:v}
  Let $1\leq p<\infty$. We define $V^p$ as the normed space of all
  functions $v:\R\to L^2$ such that $v(\infty):=\lim_{t\to
    \infty}v(t)=0$ and $v(-\infty)$ exists and for which the norm
  \begin{equation}\label{eq:hom_norm_v}
    \|v\|_{V^p}:=\sup_{\{t_k\}_{k=0}^K \in \mathcal{Z}}
    \left(\sum_{k=1}^{K}
      \|v(t_{k})-v(t_{k-1})\|_{L^2}^p\right)^{\frac{1}{p}}
  \end{equation}
  is finite.  Likewise, let $V^p_-$ denote the normed space of all
  functions $v:\R\to L^2$ such that $v(-\infty)=0$, $v(\infty)$
  exists, and $\|v\|_{V^p} < \infty$, endowed with the norm
  \eqref{eq:hom_norm_v}.
\end{d1}

\begin{p1}\label{prop:v}Let $1\leq p<q <\infty$.
  \begin{enumerate}
  \item\label{it:v_limits} Let $v:\R\to L^2$ be such that
$$\|v\|_{V_0^p}:=\sup_{\{t_k\}_{k=0}^K \in \mathcal{Z}_0}
\left(\sum_{k=1}^{K}
  \|v(t_{k})-v(t_{k-1})\|_{L^2}^p\right)^{\frac{1}{p}}
$$
is finite. Then, it follows that $v(t_0^+):=\lim_{t \downarrow
  t_0}v(t)$ exists for all $t_0\in [-\infty,\infty)$ and
$v(t_0^-):=\lim_{t \uparrow t_0}v(t)$ exists for all $t_0\in
(-\infty,\infty]$ and moreover,
$$
\|v\|_{V^p}=\|v\|_{V_0^p}.
$$
\item\label{it:v_spaces} We define the closed subspace $V^p_{rc}$
  ($V^p_{-,rc}$) of all right-continuous $V^p$ functions ($V^p_-$
  functions).  The spaces $V^p$, $V^p_{rc}$, $V^p_-$ and $V^p_{-,rc}$
  are Banach spaces.
\item\label{it:v_emb1} The embedding $U^p\subset V_{-,rc}^p$ is
  continuous.
\item\label{it:v_emb2} The embeddings $V^p\subset V^q$ and
  $V^p_-\subset V^q_-$ are continuous.
\item The embedding $V^p_{-,rc} \subset U^q$ is continuous.
\end{enumerate}
\end{p1}

\begin{t1}\label{thm:duality}
  Let $1<p<\infty$.We have
  \begin{equation*}
    (U^p)^\ast=V^{p'}
  \end{equation*}
  in the sense that there is a bilinear form $B$ such that the mapping
  \begin{equation}\label{eq:duality}
    T: V^{p'} \to (U^p)^\ast, \; T(v):=B(\cdot,v)
  \end{equation}
  is an isometric isomorphism.
\end{t1}

\begin{p1}\label{prop:c1} Let $1<p<\infty$, $u\in V^1_-$ be absolutely 
  continuous on compact intervals and $v\in V^{p'}$. Then,
  \begin{equation}\label{eq:c1}
    B(u,v)=-\int_{-\infty}^\infty \lb u'(t),v(t)\rb dt.
  \end{equation}
\end{p1}

\begin{p1}\label{cor:mod_est}
  We have, for $M = 2^k$, $k \in \Z$,
  \begin{align}\label{eq:mod_est1}
    \|Q^\pm_M u\|_{L^2(\R^n)} &\ls M^{-\frac{1}{2}} \|u\|_{V^2_\pm}\\
    \label{eq:mod_est2}
    \|Q^\pm_{\geq M} u\|_{L^2(\R^n)} &\ls M^{-\frac{1}{2}} \|u\|_{V^2_\pm}\\
    \label{eq:mod_est3}
    \|Q^\pm_{< M} u\|_{V^p_\pm} \ls \|u\|_{V^p_\pm}\ ,\quad & \|Q^\pm_{\geq M} u\|_{V^p_\pm}
    \ls \|u\|_{V^p_\pm}\\
    \label{eq:mod_est4}
    \|Q^\pm_{< M} u\|_{U^p_\pm} \ls \|u\|_{U^p_\pm}\ ,\quad &\|Q^\pm_{\geq M}
    u\|_{U^p_\pm} \ls \|u\|_{U^p_\pm}
  \end{align}
\end{p1}

\begin{p1}\label{prop:ext_free_est}
  Let $$T_0:L^2\times \cdots \times L^2\to L^1_{loc}(\R^n;\C)$$ be a
  $m$-linear operator.
  Assume that for some $1\leq
    p,q\leq\infty$
    \begin{equation*}
      \|T_0(e^{ \pm_1 \cdot \jb{D}}\phi_1,\ldots,e^{
        \pm_n \cdot \jb{D}}\phi_m)\|_{L^p_t(\R;L^q_{x,y}(\R^n))}
      \ls \prod_{i=1}^m \|\phi_i\|_{L^2}.
    \end{equation*}
    Then, there exists $T:U^p_{\pm_1}\times \cdots \times U^p_{\pm_m} \to
    L^p_t(\R;L^q_{x,y}(\R^n))$ satisfying
    \begin{equation*}
      \|T(u_1,\ldots,u_m)\|_{L^p_t(\R;L^q_{x,y}(\R^n))}\ls \prod_{i=1}^m \|u_i\|_{U_{\pm_i}^p},
    \end{equation*}
    such that $
    T(u_1,\ldots,u_m)(t)(x,y)=T_0(u_1(t),\ldots,u_m(t))(x,y) $ a.e..
\end{p1}

\begin{p1}\label{prop:interpolation} Let $q>1$,
  $E$ be a Banach space and $T:U_\pm^q\to E$ be a bounded, linear
  operator with $\|Tu\|_E \leq C_q \|u\|_{U^q_\pm}$ for all $u \in
  U^q_\pm$.  In addition, assume that for some $1\leq p <q$ there exists
  $C_p\in (0,C_q]$ such that the estimate $\|Tu\|_E \leq C_p
  \|u\|_{U^p_\pm}$ holds true for all $u \in U^p_\pm$.  Then, $T$
  satisfies the estimate
  \begin{equation*}
    \|Tu\|_E \lesssim (1+ \ln\frac{C_q}{C_p})\|u\|_{V^p_\pm},
    \quad u \in V^p_{-,rc,\pm}.
  \end{equation*}
\end{p1}
\bigskip

\subsection{Function spaces}

We rewrite equation \eqref{nlkg} (in the scalar case $K=1$, $m_1 = 1$) as a first order system, which we can more comfortably taylor our function spaces to. To this end, we note that
\[ \Box + 1 = (\jb{D} + i\partial_t)(\jb{D} - i\partial_t). \]
Hence, given a sufficiently regular function $u$ that satisfies
\[ (\Box + 1)u = F, \qquad u(0) = f,\;\; u'(0)=g \]
we define \begin{equation} \label{upm} u^\pm = \frac{\jb{D} \mp i\partial_t}{2\jb{D}}u. \end{equation}
Then the $u^\pm$ solve\footnote{the expression $\frac{F}{\jb{D}}$ corresponds to the gain of one full derivative of the linear Klein-Gordon equation}
\begin{equation} \label{eqn} (\jb{D} \pm i\partial_t)u^\pm = \frac{F}{2\jb{D}}, \qquad u^{\pm}(0) = \frac{1}{2}\left(f \mp i\frac{g}{\jb{D}}\right). \end{equation}
Since we have the identity $u^+  + u^- = u$, we may reconstruct $u$ from this system, and we will in our estimates work exclusively on \eqref{upm} and \eqref{eqn}. 

With this construction in mind, we define the function spaces which we are going to use.

\begin{d1} We define the closed spaces $X^s_\pm \subset C(\R, H^s$) by the closure of $C(\R, H^s) \cap U^2$ with respect to the norm
\[ \|u\|_{X^s_{\pm}} = \left(\sum_N N^{2s} \|P_N u^\pm\|_{U^2_{\pm}}^2 \right)^{\frac{1}{2}} \]
 \[ \text{where } \|f\|_{U^2_{\pm}} = \|e^{\mp it\jb{D}}f\|_{U^2(\R, L^2)}. \]
We also define by $Y^s$ the corresponding space where $U^2$ is replaced by $V^2 = V^2_{-, rc}$.

Furthermore, we define
\[ X^s = X^s_+ \times X^s_-, \qquad Y^s = Y^s_+ \times Y^s_-. \]

\end{d1}

With these definitions, we have
\[ X^s \subset Y^s. \]
We also define the restricted space $X^s([0, \infty))$
\[ X^s([0, \infty)) = \Big\{ u \in C([0, \infty), H^s)\ |\; \tilde{u} = \ind{[0,\infty)}(t)u(t) \in X^s \Big\} \]
with norm \[ \|u\|_{X^s([0. \infty))} = \|\ind{[0,\infty)}u\|_{X^s}. \]
and define $Y^s([0, \infty))$ analogously. They are again Banach spaces.

The strategy of the proof follows the standard approach using the contraction mapping principle. We briefly outline the procedure below.

By the equivalent formulation as a system \eqref{eqn}, a solution of \eqref{nlkg} is equivalent to
\begin{equation}\begin{aligned}\label{system} (\jb{D} \pm i\partial_t)u^\pm &= \frac{1}{2\jb{D}}	N(u^+ + u^-) \\ u^\pm(0) &= u_0^\pm \end{aligned}\end{equation}
where $u_0^\pm = \frac{1}{2}\left(f \mp i\frac{g}{\jb{D}}\right) \in B_{\epsilon}(0) \subset H^{s}$.

Hence, by a solution of the above equation, we will mean $(u^+, u^-) \in X^s([0, \infty))$ which on $[0, \infty)$ solve the operator equation 
\begin{equation}
u^\pm(t) = e^{\pm it\jb{D}}u_0^\pm \mp i I^\pm(u)
\end{equation}
where $u = u^+ + u^-$ and \begin{equation}\label{duhamel}
I^\pm(u) = \int_0^t e^{\pm i(t-s)\jb{D}}\frac{N(u(s))}{2\jb{D}} ds.
\end{equation}

This equation can be solved by a contraction mapping argument in $X^s$ once we have the bounds\footnote{along with a difference version of the nonlinear bound}
\[\begin{aligned} \|e^{\pm it\jb{D}}u_0^\pm\|_{X^s_\pm([0, \infty))} &\lesssim \|u_0^\pm\|_{H^s},\\ 
 \|I^\pm(u)\|_{X^s_\pm([0, \infty))} &\lesssim \|(u_+, u_-)\|_{X^s_\pm([0, \infty))}^2 \end{aligned}\]
The linear part of the estimate is straightforward\footnote{we brush over the restriction to $t \in [0, \infty)$, which does not harm this estimate}, since
\begin{equation} \label{linest}\|e^{\pm it\jb{D}}u_0^\pm\|_{X^s_\pm}^2 = \sum N^{2s} \|e^{\pm it\jb{D}}P_N(u_0^\pm)\|_{U^2_\pm}^2 = \sum N^{2s} \|P_N(u_0^\pm)\|_{U^2}^2 \sim \|u_0^\pm\|_{H^s}^2. \end{equation}

Hence, the focus of the sections to come is on the nonlinear estimate of $I^\pm$. In the next section, we derive some spacetime estimates that will be crucial for the nonlinear estimate.

\section{Bilinear and Strichartz estimates}
We will tacitly assume that $n \geq 2$. By virtue of \cref{prop:ext_free_est}, bounds in $U^p_\pm$ type spaces follow from $L^p$ bounds on free solutions $e^{\pm it\jb{D}}\phi$. For our estimates in dimension 3 or higher, we will use the key estimate below.
\begin{p1}\label{propl2} Let $n \geq 3$, let $O, M, N$ dyadic numbers and $\phi_M$, $\psi_N$ functions in $L^2(\R^n)$ localized at frequencies $M$, $N$ respectively. Define $u_M = e^{\pm_1it\jb{D}}\phi_M$, $v_N = e^{\pm_2 it\jb{D}} \psi_N$. Denote $L = \min(O,M,N)$, $H = \max(O, M, N)$. Then,
\begin{equation} \|P_O(u_M v_N) \|_{L^2(\R \times \R^n)} \lesssim \begin{cases} H^\frac{1}{2} L^\frac{n-2}{2}\|\phi_M\|_{L^2}\|\psi_N\|_{L^2} & \text{ if } M \sim N \\ L^{\frac{n-1}{2}}\|\phi_M\|_{L^2}\|\psi_N\|_{L^2} & \text{ otherwise} \end{cases}\end{equation}
\end{p1}
\begin{proof}see \cref{secbilin}.\end{proof}

We also state the Strichartz estimates available for the Klein-Gordon equation. These will mainly be used when $n=2$ as the above more powerful bilinear refinement is not available in that case.
The estimates come in two main flavors, depending on whether one chooses to use the radial curvature of the characteristic hypersurface. The condition $r < \infty$ serves to exclude inconvient endpoint cases, which we will not need in what follows.
\begin{p1}[KG-type Strichartz]\label{kgstrichartz}Let $2 \leq r < \infty$, $\frac{2}{q}+\frac{n}{r} = \frac{n}{2}, $ and $l = \frac{1}{q}-\frac{1}{r}+\frac{1}{2}$. Then
\begin{equation} \label{ss}\|e^{it\langle D \rangle}u_0\|_{L^q_t L^r_x} \lesssim \|\langle D \rangle^l u_0\|_{L^2}. \end{equation}
\end{p1}
\begin{proof} see \cite{MR2424390}. \end{proof}
\begin{p1}[Wave-type Strichartz] Let $2 \leq r < \infty$, $\frac{2}{q}+\frac{n-1}{r} = \frac{n-1}{2}, $ and $l = \frac{1}{q}-\frac{1}{r}+\frac{1}{2}$. Then
\begin{equation} \label{sw} \|e^{it\langle D \rangle}u_0\|_{L^q_t L^r_x} \lesssim \|\langle D \rangle^l u_0\|_{L^2}. \end{equation}
\end{p1}
\cref{propl2} implies the following bilinear refinement in $L^4$ which will be useful in controlling the worst interactions:
\begin{p1}[$L^4$ estimate]
Let $n \geq 3$ and for $M \lesssim N$, let $\phi_{N,M}$ be supported in a ball of radius $M$ located at frequency $N$. Then
\begin{equation} \label{swref} \|e^{it\jb{D}}\phi_{N,M}\|_{L^4} \lesssim N^\frac{1}{4}M^\frac{n-2}{4}\|\phi_{N,M}\|_{L^2}. \end{equation}
\end{p1}
\begin{r1} Using Wave Strichartz estimates and Bernstein's inequality, one would only obtain the weaker estimate
\[\begin{aligned} \|e^{it\jb{D}}\phi_{N,M}\|_{L^4} &\lesssim M^{(\frac{n(n-2)}{2(n-1)}-\frac{n}{4})} \|e^{it\jb{D}}\phi_{N,M}\|_{L^4 L^{\frac{2(n-1)}{n-2}}}\\
&\lesssim M^{(\frac{n(n-2)}{2(n-1)}-\frac{n}{4})} N^{\frac{3}{4}-\frac{n-2}{2(n-1)}}\|\phi_{N,M}\|_{L^2} \end{aligned}\]
where the loss in the high frequency $N$ is always strictly larger than $\frac{1}{4}$.\footnote{The same occurs using the Schr\"odinger estimate.}
\end{r1}
\begin{proof} We omit the indices $M, N$ and rewrite the estimate in the equivalent bilinear fashion
\[ \|e^{it\jb{D}}\phi e^{-it\jb{D}}\bar{\phi}\|_{L^2} \lesssim N^\frac{1}{2}M^\frac{n-2}{2}\|\phi\|_{L^2}^2 \]
Now the Fourier supports of $\phi$ and $\bar{\phi}$ are symmetric through the origin, and hence the sum of the supports is contained in a ball of radius $\lesssim M$ centered at the origin. We may thus insert a projector $P_M$ and it remains to estimate
\[ \|P_M(e^{it\jb{D}}\phi e^{-it\jb{D}}\bar{\phi})\|_{L^2} \lesssim N^\frac{1}{2}M^\frac{n-2}{2}\|\phi\|_{L^2}^2 \]
but this is one of the bilinear estimates in \cref{propl2}.
\end{proof}

With these building blocks, we transfer the estimates over on the corresponding $U^p_\pm$ and $V^2_\pm$ spaces using \cref{prop:ext_free_est} and \cref{prop:interpolation}.

\begin{p1}[$U^4 \rightarrow L^4$]\label{propu4} Let $n \geq 3$ and let $u_{M,N}$ have Fourier support in a ball of \mbox{radius $\sim M$} centered at frequency $N \gtrsim M$. Then\footnote{in the case where $u_{M,N} = P_N u$, the estimate of course still holds with $M=N$}
\begin{equation} \label{swrefu2}\|u_{M,N}\|_{L^4} \lesssim N^\frac{1}{4}M^\frac{n-2}{4}\|u_{M,N}\|_{U^4_\pm}. \end{equation}
\end{p1}
\begin{proof} This follows from \eqref{swref} and \cref{prop:ext_free_est}.
\end{proof}
\begin{p1}[$L^2 \rightarrow U^2 \times U^2$] \label{thmbilin} Let $n \geq 3$ and let $L$ ($H$) the lowest (highest) of the frequencies $M,N,O$. Let $u_M \in U^2_{\pm_1}$, $u_N \in U^2_{\pm_2}$. Then we have
\begin{equation}\label{bilin} \|P_O(u_M v_N) \|_{L^2} \lesssim \begin{cases} L^{\frac{n-1}{2}}\|u_M\|_{U^2_{\pm_1}}\|u_N\|_{U^2_{\pm_2}} & \text{ if } M \ll N \\ H^\frac{1}{2} L^\frac{n-2}{2}\|u_M\|_{U^4_{\pm_1}}\|u_N\|_{U^4_{\pm_2}} & \text{ if } M \sim N \end{cases}\end{equation}

Furthermore, we may take $u_M$ and $u_N$  in $V^2_{\pm}$ at the expense of a factor\footnote{of course $\max(1, \log(\cdot))$ is meant} $\log^2 \frac{H}{L}$ in the case $M \ll N$. When $M \sim N$, the same is true without an additional factor.
\end{p1}
\begin{proof} We omit the $\pm$ index in the $U^2$ spaces. The estimates in $U^2$ follow from \cref{thmbilin}. In the case $M \sim N$, we improve to $U^4_\pm$ by orthogonality (as outlined the proof of \cref{bilsel} in the appendix) and \cref{propu4}.
It remains to interpolate with $V^2$ when $M \ll N$. Define $Tv = P_O(u_M P_N v )$. Then we have by \eqref{swrefu2}, \eqref{bilin} and $U^2 \subset U^4$
\[ \|T\|_{U^4_{\pm_2} \rightarrow L^2} \lesssim (MN)^\frac{n-1}{4}\|u_M\|_{U^2_{\pm_1}}, \qquad \|T\|_{U^2_{\pm_2} \rightarrow L^2} \lesssim L^\frac{n-1}{2}\|u_M\|_{U^2_{\pm_1}} \]
where \[ (MN)^\frac{n-1}{4} = (HL)^\frac{n-1}{4}. \]
Since $\log \frac{(HL)^\frac{n-1}{4}}{L^\frac{n-1}{2}} \lesssim \log \frac{H}{L}$, interpolation using \cref{prop:interpolation} yields
\[ \|T\|_{V^2_{\pm_2} \rightarrow L^2} \lesssim L^\frac{n-1}{2} (\log \frac{H}{L}) \|u_M\|_{U^2_{\pm_1}} \]

Now we iterate the argument with $S: u \mapsto P_O(P_M u, v_N)$. This time we have, using $V^2 \subset U^4$,
\[ \|S\|_{U^4_{\pm_1} \rightarrow L^2} \lesssim C_{M,N}\|v_N\|_{V^2_{\pm_2}}, \qquad \|S\|_{U^2_{\pm_1} \rightarrow L^2} \lesssim L^\frac{n-1}{2} (\log \frac{H}{L}) \|v_N\|_{V^2_{\pm_2}}, \]
and hence, since $L^\frac{n-1}{2}\log \frac{H}{L} \gtrsim L^\frac{n-1}{2}$, as before
\[ \|S\|_{V^2_{\pm_1} \rightarrow L^2} \lesssim L^\frac{n-1}{2} (\log \frac{H}{L})^2 \|u_M\|_{V^2_{\pm_2}}. \]
\end{proof}

\section{Trilinear estimates}
In this section, we perform the estimates necessary to prove bounds for the Duhamel terms $I^\pm(u)$. The fact that we are dealing with quadratic nonlinearities in combination with the important duality between $U^2$ and $V^2$ - as induced by the bilinear form $B$ from \cref{thm:duality} - is why these estimates are trilinear in nature. To motivate the precise form of the proposition below, we compute with $f = \frac{N(u)}{2\jb{D}}$ for the Duhamel term \eqref{duhamel}
\begin{equation}\begin{aligned}\label{duality}
\|P_N I^\pm(u)\|_{U^2_\pm} &= \|e^{\mp it\jb{D}} I^\pm(u)\|_{U^2_0} = \|P_N \int_0^t e^{\mp is\jb{D}}f(s)ds\|_{U^2_0} \\
&= \sup_{\|v\|_{V^2}=1} \Big|B\left(P_N \int_0^t e^{\mp is\jb{D}}f(s)ds, v\right)\Big| \\
&= \sup_{\|v\|_{V^2}=1} \Big|\iint f(t)\ \overline{ e^{\pm it\jb{D}} P_N v(t)}\ dxdt\Big| \\
&= \sup_{\|P_N v\|_{V^2_\pm}=1} \Big|\iint f(t)\ \overline{P_N v(t)}\ dxdt\Big|
\end{aligned}
\end{equation}

It will become apparent shortly that \eqref{tribad} and \eqref{trigood} below are exactly the estimates needed to sum the high-low interactions and high-high interactions, respectively.

\begin{t1}[Trilinear estimates]\label{thmtrilin}
Let $s \geq \max(\frac{1}{2}, \frac{n-2}{2})$, assume that $\pm_i$ $(i~=~1,2,3)$ are arbitrary signs and that $H \sim H'$. Then,
\begin{equation}\label{tribad} \frac{1}{H} \Big| \sum_{L \lesssim H} \iint_0 u_L v_{H'} w_{H}dxdt \Big| \lesssim \left( \sum_{L \lesssim H} L^{2s}\|u_L\|_{V^2_{\pm_1}}^2 \right)^\frac{1}{2} \|v_{H'}\|_{V^2_{\pm_2}} \|w_H\|_{V^2_{\pm_3}}. \end{equation}
Also, we have
\begin{equation} \label{trigood}\left( \sum_{L \lesssim H} L^{-2} L^{2s} \sup_{\|w_L\|_{V^2_{\pm_3}} = 1} \Big|\iint_0 u_{H'}v_H w_L dxdt\Big|^2\right)^\frac{1}{2} \lesssim {H'}^s \|u_{H'}\|_{V^2_{\pm_1}} H^s \|v_H\|_{V^2_{\pm_2}}. \end{equation}
\end{t1}
\begin{proof}
The proof will use the following
\begin{l1}[Modulation bound]\label{modlemma} Let $\xi_1 + \xi_2 = \xi_3$. Then we have 
\begin{equation}\label{modcrude}
\jb{\xi_1} + \jb{\xi_2} - \jb{\xi_3} \gtrsim \jb{\xi_{min}}^{-1}.
\end{equation}
\end{l1}
The above lemma can be improved, but we will only need \eqref{modcrude}.

We decompose each function in a low and high modulation part, where the threshold between the two regimes is set at $\Lambda > 0$ which will be chosen immediately.
Recall that we defined
\[ Q^\pm_{>M}u = \F_{tx}^{-1}\left(\ind{\{|\tau \mp \jb{D}| > M\}}\F_{tx}u\right) \text{ and } Q^\pm_{\leq M} = 1 - Q^\pm_{>M}. \]
For the proof of \eqref{tribad}, we compose $u_L = u_L^h + u_L^l$, where $u_L^h = Q^{\pm_1}_{>\Lambda} u_N$. Similarly we decompose $v_{H'}$ and $w_H$, using instead the signs $\pm_2$ and $\pm_3$, respectively.
Then we have
\[ \int u_L^l v_{H'}^l w_H^l dx dt = \left(\F_{tx}u_L^l * \F_{tx}v_{H'}^l * \F_{tx}w_H^l\right)(0, 0) \]
to which only frequencies $\tau_1 + \tau_2 + \tau_3 = 0$, $\xi_1 + \xi_2 + \xi_3 = 0$ contribute. Since from the definition of $Q^\pm_{\leq \Lambda}$ we also have $\Big| \tau_i \mp_i \jb{\xi_i} \Big| \leq \Lambda$, we see that on the contributing set
\[ 3\Lambda \geq \Big| \sum_{i=1}^3 \left(\tau_i \mp_i \jb{\xi_i}\right)\Big| = \Big| \sum_{i=1}^3 \pm_i \jb{\xi_i}\Big| \gtrsim L^{-1}, \]
which is obvious when the three signs coincide and follows from \cref{modlemma} otherwise.
Chosing $\Lambda = C^{-1}L^{-1}$ for $C$ large enough will ensure that the above integral vanishes and hence in what follows, we always have high modulation on (at least) one factor. We will indicate high modulation on $f$ by $f^h$
and treat now \eqref{tribad} in the case where $u_L = u_L^h$. Namely, we estimate the term by
\[\begin{aligned} \text{LHS }\eqref{tribad} &\lesssim H^{-1}\sum_{L \lesssim H} L^\frac{1}{2}\|u_L\|_{V^2_{\pm_1}}\|P_L(v_{H'} w_H)\|_{L^2}\\ &\leq H^{-1} \left( \sum_{L \lesssim H} L^{2s}\|u_L\|_{V^2_{\pm_1}}^2\right)^\frac{1}{2} \left( \sum_{L \lesssim H} L^{1-2s} \|P_L(v_{H'} w_H)\|_{L^2}^2\right)^\frac{1}{2}. \end{aligned}\]
When $n=2$, we simply use $L^{1-2s} \leq 1$, orthogonality, and the $q = r = 4$ Strichartz estimate from \cref{kgstrichartz}, obtaining
\[ \sum_{L \lesssim H} L^{1-2s} \|P_L(v_{H'} w_H)\|_{L^2}^2 \lesssim \|v_{H'} w_H\|_{L^2}^2 \leq H^2 \|v_{H'}\|_{V^2_{\pm_2}}^2 \|w_H\|_{V^2_{\pm_3}}^2 \]
and the claim follows. When $n \geq 3$, we have by \eqref{bilin}
\[ H^{-1} \left( \sum_{L \lesssim H} L^{1-2s} \|P_L(v_{H'} w_H)\|_{L^2}^2\right)^\frac{1}{2} \lesssim H^{-1} \left( \sum_{L \lesssim H} L^{n-1-2s}H\right)^\frac{1}{2}\|v_{H'}\|_{V^2_{\pm_2}}\|w_H\|_{V^2_{\pm_3}} \]
and the claim follows since
 \[  \sum_{L \lesssim H} L^{n-1-2s} \lesssim 1+H^{n-1-2s} \lesssim H. \]
whenever $s \geq \frac{n-2}{2}$.

Now we investigate the easier case $v_{H'} = v_{H'}^h$ (the case $w_H = w_H^h$ is the same) again by putting the high modulation term in $L^2$. For $n=2$ we get the expression
\[\begin{aligned} H^{-1}\sum_{L \lesssim H} \|v_{H'}^h\|_{L^2}\|u_Lw_H\|_{L^2} \lesssim H^{-1}\sum_{L \lesssim H} L^\frac{1}{2}L^\frac{1}{2}	H^\frac{1}{2}\|u_L\|_{V^2_{\pm_1}}\|v_{H'}\|_{V^2_{\pm_2}}\|w_H\|_{V^2_{\pm_3}}& \\
\lesssim H^{-\frac{1}{2}}\left( \sum L^{2s} \|u_L\|_{V^2_{\pm_1}}^2 \right)^\frac{1}{2} \left( \sum_{L \lesssim H} L L^{1-2s}\right)^\frac{1}{2}\|v_{H'}\|_{V^2_{\pm_2}}\|w_H\|_{V^2_{\pm_3}}&
\end{aligned}\]
which gives the claim since \[ \sum_{L \lesssim H} L L^{1-2s} \leq \sum_{L \lesssim H} L \lesssim H. \]
In higher dimensions, we estimate
\[ H^{-1}\sum_{L \lesssim H} \|v_{H'}^h\|_{L^2}\|u_Lw_H\|_{L^2} \lesssim H^{-1}\sum_{L \lesssim H} L^\frac{1}{2}L^{\frac{n-1}{2}}\log^2 \frac{H}{L}\|u_L\|_{V^2_{\pm_1}}\|v_{H'}\|_{V^2_{\pm_2}}\|w_H\|_{V^2_{\pm_3}}. \]
After Cauchy-Schwarz with $L^s\|u_L\|_{V^2}$ and the rest, using $\log^4 \frac{H}{L} \lesssim \frac{H}{L}$,
\[ H^{-2} \sum_{L \lesssim H} L L^{n-1}L^{-2s}\log^4 \frac{H}{L} \lesssim H^{-1} \sum_{L \lesssim H}L^{n-1-2s} \lesssim 1. \]
We now turn to the proof of \eqref{trigood} and perform the same modulation decomposition as before, starting with the case $w_L = w_L^h$. Then, for $n = 2$, using the high modulation, orthogonality and finally Strichartz estimates,
\[\begin{aligned} \eqref{trigood}^2 &\leq \sum_{L \lesssim H} L^{2s-1}\|P_L(u_{H'}v_H)\|^2 \lesssim H^{2s-1} \|u_{H'}v_H\|_{L^2}^2 \lesssim H^{2s+1} \|u_{H'}\|_{V^2_{\pm_1}}^2\|v_H\|_{V^2_{\pm_2}}^2\\ &\lesssim H'^{2s}\|u_{H'}\|_{V^2_{\pm_1}}^2 H^{2s} \|v_H\|_{V^2_{\pm_2}}^2 \end{aligned}\]
and for $n \geq 3$
\[ \eqref{trigood}^2 \leq \sum_{L \lesssim H} L^{2s-1}\|P_L(u_{H'}v_H)\|^2 \lesssim \sum_{L \lesssim H} L^{2s+n-3}H \|u_{H'}\|_{V^2_{\pm_1}}^2\|v_H\|_{V^2_{\pm_2}}^2. \]
Now \[  H\sum_{L \lesssim H} L^{2s+n-3} \lesssim H^{2s+n-2}\leq H'^{2s}H^{2s} \]
and the claim follows. Finally, we treat the last case $u_{H'} = u_{H'}^h$. For $n = 2$, we get
\[ \begin{aligned} \eqref{trigood}^2 &\leq \sum_{L \lesssim H} L^{2s-1}\|v_{H}w_L\|^2\|u_{H'}\|_{V^2_{\pm_1}}^2  \lesssim \sum_{L \lesssim H} L^{2s}H \|u_{H'}\|_{V^2_{\pm_1}}^2\|v_H\|_{V^2_{\pm_2}}^2 \\ &\lesssim H'^{2s}\|u_{H'}\|_{V^2_{\pm_1}}^2 H^{2s} \|v_H\|_{V^2_{\pm_2}}^2 \end{aligned}\]
and for $n \geq 3$, we obtain
\[\begin{aligned} \eqref{trigood}^2 &\lesssim \sum_{L \lesssim H} L^{2s-2} \sup_{\|w\|_{V^2}=1} L\|u_{H'}\|_{V^2_{\pm_1}}^2\|v_H w_L\|_{L^2}^2 \\&\lesssim \sum_{L \lesssim H} L^{2s-2+n}(\log^4 \frac{H}{L}) \|u_{H'}\|_{V^2_{\pm_1}}^2\|v_H\|_{V^2_{\pm_2}}^2 \end{aligned}\]
but we can replace $\log^4 \frac{H}{L}$ by $\frac{H}{L}$ and estimate as in the last case.
\end{proof}

\section{Proof of the main theorem}

Recall that we have a nonlinearity which is a finite sum of quadratic terms in $u^+$, $u^-$ and their conjugates. For brevity, we will from now on restrict to one such term without conjugates. Since $\|\bar{v}\|_{X^s_\pm} = \|v\|_{X^s_\mp}$, the other cases follow in the same manner. The main result is
\begin{t1}\label{nonlinest} let $s \geq \max(\frac{1}{2},\frac{n-2}{2})$. For any $\pm_1$, $\pm_2$, we have
\[ I_{\pm_1 \pm_2} := I:\quad Y^s \times Y^s \rightarrow X^s, \]
where \[\begin{aligned} I((u^+,u^-),(v^+, v^-)) &= (I^+(u^{\pm_1}, v^{\pm_2}), I^-(u^{\pm_1}, v^{\pm_2})) \\
I^\pm(f,g)  &= \int_0^t e^{\pm i(t-s)\jb{D}}\frac{f g}{2\jb{D}} ds. \end{aligned}\]
In other words, for a constant $C = C(n)$,
\[ \|I(u,v)\|_{X^s} \leq C\|u\|_{Y^s}\|v\|_{Y^s}. \]
In particular, since $X^s \subset Y^s$, we also have
\[ I: X^s \times X^s \rightarrow X^s \]
and
\[ I: Y^s \times Y^s \rightarrow Y^s. \]
\end{t1}

\begin{proof}
It suffices to consider the terms $S_i$, $i = 1,2$, where
\[ S_1 = \Big\|\sum_{H} \sum_{L \ll H} I(\vec{u}_L, \vec{v}_H)\Big\|_{X^s},\qquad S_2 = \Big\|\sum_{H} \sum_{H \sim H'} I(\vec{u}_H, \vec{v}_H')\Big\|_{X^s}. \]
We treat $I^+$ only since the other component follows in the same manner, denote by $u_H$ and $v_H$ the components of $\vec{u}_H$ and $\vec{v}_H$ as selected by the signs $\pm_1$ and $\pm_2$, and begin by estimating $S_1$. By duality and \eqref{tribad} from \cref{thmtrilin},
\[\begin{aligned} \Big\|P_H \sum_{L \ll H} I^+(u_L, v_H)\Big\|_{U^2_+} &= \frac{1}{H}\sup_{\|w_{H'}\|_{V^2_-}=1} \Big| \sum_{L \ll H} \iint u_L v_H w_{H'}dxdt \Big| \\  &\lesssim \left( \sum_{L \lesssim H} L^{2s}\|u_L\|_{V^2_{\pm_1}}^2 \right)^\frac{1}{2}\|v_H\|_{V^2_{\pm_2}} \end{aligned}\]
and thus
\[ \sum_H H^{2s}\Big\|P_H \sum_{L \ll H} I^+(u_L, v_H)\Big\|_{U^2_+}^2 \lesssim \|\vec u\|_{Y^s}^2\|\vec v\|_{Y^s}^2. \]
For $S_2$, we instead estimate
\[ S_2 \leq \sum_H \sum_{H' \sim H}\| I^+(u_{H'}, v_H)\|_{X^s_+} \lesssim \sum_H \sum_{H' \sim H}\left( \sum_{L \lesssim H} L^{2s} \|P_L I^+(u_{H'}, v_H)\|_{U^2_+}^2 \right)^\frac{1}{2}. \]
Using duality again, we arrive exactly at $\sum_H \sum_{H' \sim H} \eqref{trigood}$, and using this, we get
\[ S_2 \lesssim \sum_H \sum_{H' \sim H} H'^{2s}\|u_{H'}\|_{V^2_{\pm_1}} H^{2s}\|v_H\|_{V^2_{\pm_2}} \lesssim \|\vec u\|_{Y^s}\|\vec v\|_{Y^s}. \]
\end{proof}

We now solve \eqref{nlkg} by contraction mapping techniques, e.g. we are going to construct a solution of the operator equation
\begin{equation}\label{inteqn}
u^\pm(t) = T^\pm u^\pm := e^{\pm it\jb{D}}u_0^\pm \mp i I^\pm(u)
\end{equation}
where $u = u_+ + u_-$ and \begin{equation}
I^\pm(u) = \int_0^t e^{\pm i(t-s)\jb{D}}\frac{N(u(s))}{2\jb{D}} ds.
\end{equation}

We look for the solution in the set \[ D_\delta = \{ u \in X^s([0, \infty):\ \|u\|_{X^s([0, \infty))} \leq \delta \}. \]
For $u \in D_\delta$ and initial data $u_0 = (u_0^+, u_0^-)$ of size at most $\eps = \eps(\delta) \ll \delta$, we have
\[ \|e^{\pm it\jb{D}}u_0^\pm \mp i I^\pm(u)\|_{X^s_\pm([0, \infty))} \lesssim \eps + \delta^2 \leq \delta \]
for small enough $\delta$, due to \eqref{linest} and the fact that $I^\pm(u)$ is a sum of operators for which \cref{nonlinest} holds.
Since we can factor $a^2-b^2 = a(a-b) + (a-b)b$, we also obtain
\[\begin{aligned} \|I^\pm(f)-I^\pm(g)\|_{X^s_\pm([0, \infty))} &\lesssim (\|f\|_{X^s([0, \infty))}+\|g\|_{X^s([0, \infty))})\|f-g\|_{X^s([0, \infty))} \\ &\lesssim \delta \|f-g\|_{X^s([0, \infty))} \end{aligned}\]
and hence $T$ is a contraction on $D_\delta$ when $\delta \ll 1$, which implies the existence of a unique fixed point in $D_\delta$ solving the integral equation \eqref{inteqn}.

As for scattering, by \cref{nonlinest} we have that for each $N$, \[ e^{\mp it\jb{D}}P_N I^\pm(u) \in V^2_{-,rc}\] and hence, the limit as $t \rightarrow \infty$ exists for each piece. Together with
\[ \sum_N N^{2s} \|P_N I^\pm(u)\|_{V^2_\pm}^2 \lesssim 1, \]
it follows that $\lim_{t \rightarrow \infty} e^{\mp it\jb{D}} I^\pm(u) \in H^s$.
Hence, for the solution $u = (u^+, u^-)$ we have that
\[ e^{\mp it\jb{D}} u^\pm \rightarrow u_0^\pm \mp i \lim_{t \rightarrow \infty} e^{\mp it\jb{D}} I^\pm(u) \in H^s. \]

\section{Systems of different masses}\label{diffmass}
Since the bilinear estimates easily tolerate interactions between waves with different masses when $n \geq 3$ and the case $n = 2$ relies on Strichartz estimates only, the only obstruction to carrying out the proof of the main result for a system of such type is the absence of resonances\footnote{in the framework of space-time resonances (cf. \cite{germain}), our notion describes the absence of time resonance}.
Recalling the notation $\jb{\cdot} = \sqrt{m^2+|\cdot|^2}$, we have the following

\begin{l1} Let positive masses $m_1, \ldots, m_N$ be given such that for any triple \[ (m,n,o) \in \left(\{ m_i \}_{i=1}^N\right)^3 \]
we have
\begin{equation}\label{suffsys} m + n > o \end{equation}
Then we have the modulation bound
\begin{equation}\label{modbound} \jb{\xi}_m + \jb{\eta}_n - \jb{\xi+\eta}_o \gtrsim \jb{\min(|\xi|, |\eta|, |\xi+\eta|)}^{-1}. \end{equation}
\end{l1}
\begin{r1}The condition $m + n > o$ is similar to (albeit more restrictive\footnote{the fact that we need positivity as opposed to nonvanishing of this expression seems related to the fact that our method does not take advantage of the absence of space resonances} than) the condition \[ |m_1 + m_2 - m_3| \neq 0, \] which appears in numerous places, most recently in \cite{ionescupausader}.

Of course \eqref{suffsys} is equivalent to \[ 2\min \{ m_i \} > \max \{ m_i \}. \]
\end{r1}
The statements of \cref{mainthm} follow by inspection of the main arguments if in analogy to \cref{modlemma} we have the modulation bound \eqref{modbound} by obvious adaption of the function spaces and estimates to systems. We omit the details; it remains to prove \eqref{modbound}.
\begin{proof}
By symmetry, we may assume $|\eta| \leq |\xi|$. Expanding the left hand side of  \eqref{modbound} with \[ \Lambda := \jb{\xi}_m + \jb{\eta}_n + \jb{\xi+\eta}_o \sim \jb{\xi}, \] it remains to look at the expression
\[ m^2 + n^2 - o^2 + 2\jb{\xi}_m\jb{\eta}_n - 2\xi\cdot\eta. \]
If $\xi \cdot \eta \leq 0$ then, since $m + n - o > 0$,
\[\begin{aligned} m^2 + n^2 - o^2 + 2\jb{\xi}_m\jb{\eta}_n &\geq m^2 + n^2 - o^2 + 2\max(mn, \jb{\xi}_m\jb{\eta}_n)\\ &\gtrsim \jb{\xi}\jb{\eta}, \end{aligned} \]
we have \[ \jb{\xi}_m + \jb{\eta}_n - \jb{\xi+\eta}_o \gtrsim \frac{\jb{\xi}\jb{\eta}}{\Lambda} \gtrsim \jb{\eta} \]
which implies the claim regardless of whether $|\xi+\eta|$ or $|\eta|$ is the smallest number. Hence it remains to deal with the case where $\xi \cdot \eta  > 0$, in which $|\eta|$ is comparable to the minimum frequency, and we replace $\xi \cdot \eta$ by $|\xi||\eta|$ to deal directly with the worst case\footnote{if one does not do this, one sees that the worst case happens for interactions along a line and the general case is much better, but we ignore this here}.
With some hindsight, we rewrite the resulting expression as 
\[ (m+n)^2 - o^2 - 2\eps mn+ 2\left(\jb{\xi}_m\jb{\eta}_n - |\xi||\eta| - (1-\eps)mn\right) \]
where we chose $\eps \ll 1$ such that \[ (m+n)^2-o^2 -2\eps mn > 0, \]
and we now prove that \[ \jb{\xi}_m\jb{\eta}_n - |\xi||\eta| - (1-\eps)mn \]
is nonnegative and has the correct growth. We rewrite as
\[ \jb{\xi}_m\jb{\eta}_n - |\xi||\eta| - (1-\eps)mn = \frac{m^2n^2 + n^2|\xi|^2 + m^2|\eta|^2 - (1-\eps)mn\jb{\xi}_m\jb{\eta}_n -(1-\eps)mn|\eta||\xi|}{\jb{\xi}_m\jb{\eta}_n + |\xi||\eta|}\]
and estimate the nominator using $ab \leq \frac{1}{2}(a^2+b^2)$ from below by
\[\begin{aligned}
&m^2n^2 + n^2|\xi|^2 + m^2|\eta|^2 - \frac{(1-\eps)}{2}\left(n^2\jb{\xi}_m^2 + m^2\jb{\eta}_n^2 + n^2|\xi|^2 + m^2|\eta|^2\right)\\
=\ &m^2n^2 + n^2|\xi|^2 + m^2|\eta|^2 - (1-\eps)\left(n^2m^2 + n^2|\xi|^2 + m^2|\eta|^2\right)\\
\gtrsim \ &\eps\jb{\xi}^2.
\end{aligned}\]

Hence, we have bounded \[ \jb{\xi}_m + \jb{\eta}_n - \jb{\xi+\eta}_o \gtrsim \frac{\jb{\xi}^2}{\jb{\xi}\jb{\eta}\Lambda} \gtrsim \jb{\eta}^{-1} \]
as claimed.
\end{proof}

\appendix

\section{Proof of the bilinear estimates}\label{secbilin}
In this sections we prove the bilinear estimates, which are essentially identical to those of the free wave equation. We assume $n \geq 3$ and also allow for different masses to be able to treat more general systems. For this, we define
\begin{d1} \[ \jb{\cdot}_m = \sqrt{m^2 + |\cdot|^2}. \] \end{d1}
Throughout this section, we denote
\[ \|\cdot\|_{L^2} = \|\cdot\|_{L^2_{tx}} = \|\cdot\|_{L^2(\R \times \R^n)}, \quad \|\cdot\|_{L^2_x} = \|\cdot\|_{L^2(\R^n)}. \]
\begin{p1} Let $n \geq 3$, let $O, M, N \geq 1$ dyadic numbers and $\phi_M$, $\psi_N$ functions in $L^2_x$ localized at frequencies $M$, $N$ respectively. Define $u_M = e^{\pm_1it\jb{D}_{m_1}}\phi_M$, $v_N = e^{\pm_2 it\jb{D}_{m_2}} \psi_N$. Denote $L = \min(O,M,N)$, $H = \max(O, M, N)$. Then,
\begin{equation} \|P_{\leq O}(u_M v_N) \|_{L^2} \lesssim \begin{cases} H^\frac{1}{2} L^\frac{n-2}{2}\|\phi_M\|_{L^2_x}\|\psi_N\|_{L^2_x} & \text{ if } M \sim N \\ L^\frac{n-1}{2}\|\phi_M\|_{L^2_x}\|\psi_N\|_{L^2_x} & \text{ otherwise} \end{cases}\end{equation}
\end{p1}
\begin{r1} One could easily improve the constant in the first case to $L^{n-1}$ if the signs $\pm_1$ and $\pm_2$ coincide, but we do not pursue this here.
\end{r1}
\begin{proof} The statements follow from \cref{bilsel} below upon approximation of \[ \delta(\tau\pm\jb{\xi}_{m}) \] by $\eps^{-1} \ind{|\tau-\jb{\xi}_m| \leq \eps}$ when $L \gg 1$. It remains to deal with the part where $L \lesssim 1 \ll H$. This is a routine exercise due to the fact that in that case, uniformly transversal hypersurfaces interact on a region of diameter $L$. We omit the details.
\end{proof}

We will generally follow the strategy carried through for $n=3$ in \cite{MR2439535} for the wave equation, relying solely on estimates of intersections of thickened spheres. At high frequencies, the characteristic surface resembles the cone, and we stay in this regime due to the condition $L \gg 1$.

\begin{l1}\label{slemma} Let $n \geq 3$, $0 < \delta, \Delta \ll 1 \lesssim \min(r, R, L)$ and define
\[ S_\delta(r) = \{ \xi \in \R^n:\ r-\delta \leq |\xi| \leq r+\delta \}. \]
Then, for $|\xi_0| \gtrsim \max(r,R)$, and denoting by $T(\xi, L)$ the tube of radius $L$ in the direction of $\xi_0$, we have
\[ \left| T_L(\xi_0) \cap S_\delta(r) \cap \left( \xi_0 + S_\Delta(R) \right)\right| \lesssim \frac{\min(r, R, L)^{n-3}r R \delta\Delta}{|\xi_0|}. \]
\end{l1}
\begin{r1} The statement is symmetric in $(r,\delta)$ and $(R,\Delta)$. \end{r1}
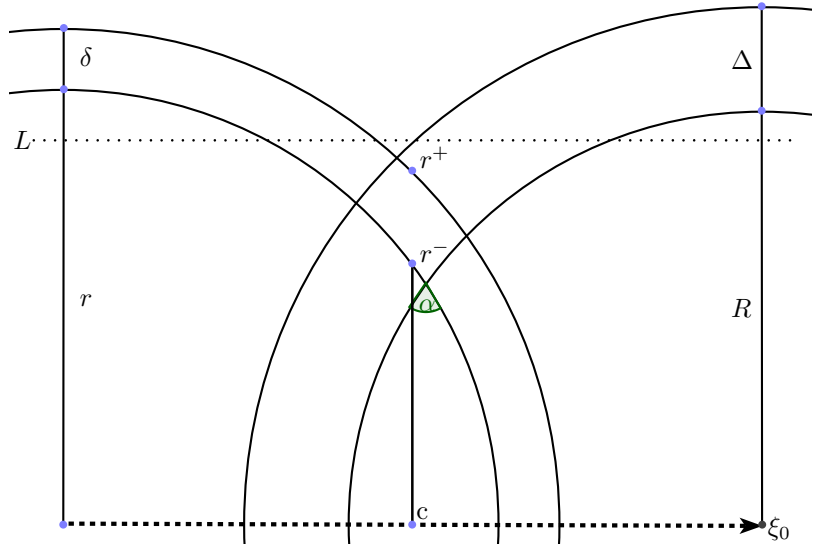
\begin{figure}[htb]
\centering
\newrgbcolor{xdxdff}{0.49 0.49 1}
\newrgbcolor{qqwuqq}{0 0.39 0}
\psset{xunit=0.5cm,yunit=0.5cm,algebraic=true,dotstyle=o,dotsize=3pt 0,linewidth=0.8pt,arrowsize=3pt 2,arrowinset=0.25}
\begin{pspicture*}(-20,-0.5)(1.3,14.5)
\pscircle(0,0){6.9}
\pscircle(0,0){5.51}
\pscircle(-18.58,0){5.8}
\pscircle(-18.58,0){6.61}
\psline(-18.58,0)(-18.57,11.59)
\psline(-18.57,11.59)(-18.57,13.22)
\psline(0,0)(-0.01,11.02)
\psline(-0.01,11.02)(-0.01,13.8)
\psline[linewidth=1.6pt,linestyle=dashed,dash=2pt 2pt]{->}(-18.44,0.03)(0.01,-0.03)
\psline[linewidth=1.0pt,linestyle=dotted]{}(-19.4,10.23)(1.01,10.23)
\rput[bl](-19.9,10){$L$}
%\pscustom[linecolor=qqwuqq,fillcolor=qqwuqq,fillstyle=solid,opacity=0.1]{\parametricplot{-0.7931225988514464}{0.8810445383627682}{0.85*cos(t)+-10.78|0.85*sin(t)+8.59}\lineto(-10.78,8.59)\closepath}
\pscustom[linecolor=qqwuqq,fillcolor=qqwuqq,fillstyle=solid,opacity=0.1]{\parametricplot{-2.1690175527527824}{-1.0449940853192858}{0.78*cos(t)+-8.94|0.78*sin(t)+6.44}\lineto(-8.94,6.44)\closepath}
\psline(-9.3,6.95)(-9.3,0)
\psline(-9.3,6.95)(-9.3,0)
\psdots[dotstyle=*,linecolor=darkgray](0,0)
\rput[bl](0.15,-0.40){$\xi_0$}
\psdots[dotstyle=*,linecolor=xdxdff](-18.58,0)
\psdots[dotstyle=*,linecolor=xdxdff](-18.57,11.59)
\rput[bl](-18.14,5.81){$r$}
\psdots[dotstyle=*,linecolor=xdxdff](-18.57,13.22)
\rput[bl](-18.14,12.22){$\delta$}
\psdots[dotstyle=*,linecolor=xdxdff](-0.01,11.02)
\rput[bl](-0.82,5.53){$R$}
\psdots[dotstyle=*,linecolor=xdxdff](-0.01,13.8)
\rput[bl](-0.82,12.12){$\Delta$}
\rput[bl](-9.13,5.70){\qqwuqq{$\alpha$}}
\psdots[dotstyle=*,linecolor=xdxdff](-9.3,6.95)
\rput[bl](-9.10,7.00){$r^-$}
\psdots[dotstyle=*,linecolor=xdxdff](-9.3,0)
\rput[bl](-9.19,0.16){c}
\psdots[dotstyle=*,linecolor=xdxdff](-9.3,9.42)
\rput[bl](-9.10,9.5){$r^+$}
\end{pspicture*}
\caption{The worst case in the proof of \cref{slemma} when $r \leq R$ in the case $n \geq 3$. Near height $L$, we have $\alpha \sim \frac{L}{r}$ and hence the intersection has volume $ \alpha^{-1}\delta \Delta L^{n-2} \sim rL^{n-3}\delta\Delta$. We remark that in two dimensions, the critical intersection occurs as the circles touch tangentially, resulting in worse estimates.}
\label{fig:circ}
\end{figure}

\begin{proof} Denote by $A$ the above intersection. We may assume $L \lesssim r$, $|\xi_0|\sim \max(r,R)$ and \[ \xi_0 = (|\xi_0|, 0, \ldots, 0). \] Hence $
\xi \in S_\delta(r) \cap (\xi_0 + S_\Delta(R))$ if and only if
\[ (r-\delta)^2 < (\xi^1)^2 + |\xi'|^2 < (r+\delta)^2 \]
and
\[ (R-\Delta)^2 < (\xi^1-\xi_0^1)^2 + |\xi'|^2 < (R+\Delta)^2. \]
Subtracting these inequalities, we find that
\[ (r-\delta)^2-(R+\Delta)^2 < (\xi^1)^2 - (\xi^1-|\xi_0|)^2 < (r+\delta)^2 - (R-\Delta)^2 \]
and hence that $\xi^1 \in (a,b)$, where
\[\begin{aligned}  a &= \frac{1}{2|\xi_0|}\left(|\xi_0|^2 + r^2-R^2 + \Delta^2-\delta^2 - 2(\delta r + \Delta R) \right) \\
b &= a + \frac{2}{|\xi_0|}(r\delta+R\Delta). \end{aligned}
\]
In particular,
\[ b-a \sim \frac{\max(r\delta,R\delta)}{|\xi_0|} \]
and hence it suffices to show, for $c \in (a,b)$,
\[ v(c) := H^{n-1}(A \cap \{ \xi^1 = c \}) \lesssim \min(r,R,L)^{n-3}\min(r\delta,R\Delta). \]
We now define the upper and lower radius over the slice $\{ \xi^1 = c \}$ of $\xi_0+S_\Delta(R)$ and $S_\delta(r)$ respectively by
\[\begin{aligned}
R^\pm(c) &= \max(0, \sqrt{(R \pm \Delta)^2 - (c-|\xi_0|)^2}) \\
r^\pm(c) &= \max(0, \sqrt{(r \pm \Delta)^2 - c^2})
\end{aligned}
\]
Ignoring for a second the intersection and only looking at $\xi_0 + S_\Delta(R)$, polar coordinates when $R^- \neq 0$ give
\[ v(c) \lesssim R^+(c)^{n-2} (R^+(c)-R^-(c)) = R^+(c)^{n-2} \frac{R^+(c)^2-R^-(c)^2}{R^+(c)+R^-(c)} \sim R^+(c)^{n-3}R\Delta, \]
whereas in the case $R^- = 0$ we get the better estimate
\[ v(c) = R^+(c)^{n-1} \lesssim (R\Delta)^\frac{n-1}{2} = (R\Delta) (R\Delta)^\frac{n-3}{2} \]
since $R+\Delta > c > R-\Delta$.

Of course the same can be done for $S_\delta(r)$, resulting in
\[ v(c) \lesssim r^+(c)^{n-3}r\delta. \]
Due to the tube $T_L(\xi_0) = T_L( (1,0, \ldots, 0))$ we also know that $\max(r^+, R^+) \leq L$, furthermore, of course, $\max(r^+, R^+) \lesssim \min(r,R)$ on $A$. In combination,
\[ v(c) \lesssim \min(r,R,L)^{n-3}\min(r\delta, R\Delta) \] as claimed, and we can estimate
\[ |A| \leq \int_a^b v(c) dc \leq (b-a)\min(r,R,L)^{n-3}\min(r\delta,R\Delta) \sim \frac{r\delta R\Delta \min(r,R,L)^{n-3}}{|\xi_0|}. \]
\end{proof}

\begin{d1} For $M, m, \eps > 0$, we denote \[ K_{M, \eps}^{m, \pm} = \left\{ (\tau, \xi) \in \R \times \R^n:\, |\xi| \sim M, \, |\tau-\jb{\xi}_m| \leq \epsilon \right\} \]
\end{d1}
\begin{p1}\label{bilsel} Let $n \geq 3$, let $O,M,N > 1$ dyadic numbers, denote \[ H = \max(O,M,N),\qquad L = \min(O,M,N)\] and let \[ \supp(u) \subseteq K_{M, \eps_1}^{m_1, \pm_1},\; \; \supp(v) \subseteq K_{N, \eps_2}^{m_2,\pm_2}. \]
Then we have
\begin{equation} \|P_{\leq O}(uv)\|_{L^2} \lesssim \begin{cases} (\eps_1 \eps_2)^\frac{1}{2}H^\frac{1}{2} L^\frac{n-2}{2}\|u\|_{L^2}\|v\|_{L^2} & \text{ if } M \sim N \\ (\eps_1 \eps_2)^\frac{1}{2}L^\frac{n-1}{2}\|u\|_{L^2}\|v\|_{L^2} & 1 \ll M \ll N \end{cases}\end{equation}
\end{p1}
\begin{r1} Note that for technical reasons, we do not treat $1 \sim M \ll N$ here.
\end{r1}
\begin{proof} If $H \lesssim 1$, losing derivatives does not matter, and the statement follows easily. Hence, in what follows, we may assume $H \gg 1$.

\subsection*{case 1} We begin with the case where $M \sim N \sim H \gg L \sim O$. Decomposing the spatial frequency supports of $u$ and $v$ in balls of radius $L$, we note that it suffices to prove the estimate
\begin{equation}\label{bilinred1}
\|P_{\leq L}(u^B v^{B'})\|_{L^2} \lesssim (\eps_1 \eps_2)^\frac{1}{2}H^\frac{1}{2} L^\frac{n-2}{2}\|u\|_{L^2}\|v\|_{L^2},
\end{equation}
where $u^B$ and $v^{B'}$ are supported in balls $B, B'$ of radius $L$ located at frequency $H$. Indeed, denote by $\mathcal{B}$ a reasonable covering of $\{ |\xi| \sim H \}$ with such balls. Then we can estimate
\[ \|P_{\leq O}(u v)\|_{L^2} \lesssim \sum_{B,B' \in \mathcal{B}} \|P_{\leq O}(u^B v^{B'})\|_{L^2} \sim \sum_{B \sim B'}  \|u^B v^{B'}\|_{L^2} \]
where $B \sim B'$ if and only if $(B+B') \cap B(0, L) \neq \emptyset$. Since for fixed $B$ there are only finitely many $B'$ with $B \sim B'$, we can further estimate
\[ \sum_{B \sim B'}  \|u^B v^{B'}\|_{L^2} \lesssim (\eps_1 \eps_2)^\frac{1}{2}H^\frac{1}{2} L^\frac{n-2}{2}\sum_{B \sim B'} \|u^B\|_{L^2}\|v^B\|_{L^2} \lesssim (\eps_1 \eps_2)^\frac{1}{2}H^\frac{1}{2} L^\frac{n-2}{2} \|u\|_{L^2}\|v\|_{L^2}. \]
Hence, we only need to prove \eqref{bilinred1}, but we may even reduce to the case $u^B = v^{B'}$ since
\[ \|u^B v^{B'}\|_{L^2}^2 \leq \|(u^B)^2\|_{L^2}\|(v^{B'})^2\|_{L^2}. \]
Furthermore, we may assume $\pm_1 = +$, $m = 1$.
We will need the following well-known
\begin{l1} Let $\supp \F_{tx} u \subseteq A$, $\supp F_{tx} v \subseteq B$. Then
\[ \|uv\|_{L^2} \leq \left( \sup_{\tau, \xi} |A \cap ((\tau, \xi)-B)|\right)^\frac{1}{2} \|u\|_{L^2} \|v\|_{L^2}. \]
\end{l1}
\begin{proof}
We denote $\F_{tx} \cdot = \tilde{\cdot}$, $\zeta = (\tau, \xi)$, $\zeta' = (\tau', \xi')$ and estimate
\[\begin{aligned} \|uv\|_{L^2}^2 &= \|\tilde{u} * \tilde{v}\|_{L^2}^2 = \|\ind{A}\tilde{u} * \ind{B}\tilde{v}\|_{L^2}^2 \\
&= \int \left( \int (\ind{B}\tilde{v})(\zeta-\zeta')(\ind{A}\tilde{u})(\zeta') d\zeta' \right)^2 d\zeta \\
&= \int \left( \int (\ind{A \cap (\zeta - B)}(\zeta') \tilde{v}(\zeta-\zeta')\tilde{u}(\zeta') d\zeta' \right)^2 d\zeta \\
&\leq \int \left( \int \ind{A \cap (\zeta - B)}(\zeta') d\zeta'\right)\left(\int |\tilde{v}(\zeta-\zeta')\tilde{u}(\zeta')|^2 d\zeta' \right) d\zeta \\
&\leq \left( \sup_\zeta |A \cap (\zeta - B)| \right) \|u\|_{L^2}^2\|v\|_{L^2}^2.
\end{aligned}\]
\end{proof}
Applying the lemma to the present situation, where $A = B$ is a ball of radius $L$ located at frequency $H$, we see that the constant in the estimate is $\sqrt{|E|}$, where
\[\begin{aligned} E = \Big\{ (\tau, \xi):\ \xi \in B,\ \xi_0-\xi \in B,\ \tau &= \jb{\xi} + O(\eps_1),\\ \tau_0-\tau&=\jb{\xi_0-\xi} + O(\eps_1)\Big\} \end{aligned}\]
uniformly in $(\tau_0, \xi_0)$. We denote \[ E(\tau) = \left\{ \xi:\ (\tau, \xi) \in E \right\} \]
and note that since $H \gg 1$, we have $\jb{\xi} \sim |\xi|$ and hence $E(\tau) = \emptyset$ unless \[ \tau = |\text{center(B)}| + O(L). \]  Thus, we have
\[ |E| \lesssim L\sup_{\tau} |E(\tau)|. \]
For $\xi \in E(\tau)$ we have $\xi_0 = (\xi_0-\xi)+\xi \in B+B \subset \{ |\eta| \sim H \}$. Now we note that
\[ |\tau-\jb{\xi}| \leq \eps_1 \iff |\xi| \in \left[ \sqrt{(\tau-\eps_1)^2-1}, \sqrt{(\tau+\eps_1)^2-1} \right]. \]
This interval has length comparable to $\eps_1$ and contains $\sqrt{\tau^2-1}$, hence it follows that
\[ |\tau-\jb{\xi}| \leq \eps_1 \Rightarrow \xi \in S_{C\eps_1}(\sqrt{\tau^2-1}). \]
In the same way,
\[ |(\tau_0-\tau)-\jb{\xi_0-\xi}| \leq \eps_1 \Rightarrow \xi_0-\xi \in S_{C\eps_1}(\sqrt{(\tau_0-\tau)^2-1}) \]
so that \[ E(\tau) \subset S_{C\eps_1}(\sqrt{\tau^2-1}) \cap \left(\xi_0 + S_{C\eps_1}(\sqrt{(\tau_0-\tau)^2-1})\right). \]
Remembering the additional restriction that the intersection happens in the ball $B$ of radius $L$, and that $\xi_0 \in 2B$, we may intersect this last set with the tube of radius $L$ along $\xi_0$. This puts us right in the situation of \cref{slemma} about intersections of thin shells, and noting that $\sqrt{\tau^2-1} \sim H \sim \sqrt{(\tau_0-\tau)^2-1}$ together with $|\xi_0| \sim H$ gives
\[ |E| \lesssim L\sup_\tau |E(\tau)| \lesssim L \frac{H^2 L^{n-3}\eps_1\eps_2}{|\xi_0|} \sim L^{n-2}H \eps_1\eps_2. \]

\subsection*{case 2} Now, without loss of generality, $L \sim M$, $H \sim N \sim O$ and $\pm_1 = +$. We may replace the projector $P_{\leq O}$ by a projector on an annulus $P_O$ (see \cite{MR2439535}, 4.3.3). Again, we want to estimate $|E|$, where
\[ E = K_{L, \eps_1}^{m_1, +} \cap \left( (\tau_0, \xi_0) - K_{N, \eps_2}^{m_2, \pm_2}\right) \]
and we have $|\xi_0| \sim H$ due to the projector $P_O$. Going through the same procedure as before, we obtain
\[\begin{aligned} E = \Big\{ (\tau, \xi):\ |\xi| \sim L, |\xi_0-\xi| \sim H,\ \tau &= \jb{\xi}_{m_1} + O(\eps_1),\\ \tau_0-\tau&=\pm \jb{\xi_0-\xi}_{m_2} + O(\eps_2)\Big\}\end{aligned} \]
and, recalling that $L \gg 1$, \[ E \subset S_{C\eps_1}(\sqrt{\tau^2-m_1^2}) \cap \left(\xi_0 + S_{C\eps_2}(\sqrt{(\tau_0-\tau)^2-m_2^2})\right). \]
Now we have 
\[ \sqrt{\tau^2-m_1^2} \sim L,\quad \sqrt{(\tau_0-\tau)^2-m_2^2} \sim H, \quad |\xi_0| \sim H\] and thus
\[ |E| \lesssim L\sup_\tau |E(\tau)| \lesssim L\frac{HL L^{n-3}\eps_1 \eps_2}{H} = L^{n-1}\eps_1 \eps_2 \]
which gives the constant $L^\frac{n-1}{2}(\eps_1 \eps_2)^\frac{1}{2}$ as claimed.

\end{proof}

\bibliographystyle{amsalpha}
\bibliography{KG2D}
\end{document}